\nonstopmode
\documentclass[10pt]{amsart}

\usepackage{amssymb, amsmath, euscript}

\input epsf

\usepackage{amscd}
\usepackage{graphicx}

\begin{document}


\newtheorem{theorem}{Theorem}[section]
\newtheorem*{theorem*}{Theorem}
\newtheorem*{observation*}{Observation}
\newtheorem*{lemma*}{Lemma}
\newtheorem{thm}[theorem]{Theorem}
\newtheorem{proposition}[theorem]{Proposition}
\newtheorem{prop}[theorem]{Proposition}
\newtheorem{lemma}[theorem]{Lemma}
\newtheorem{lem}[theorem]{Lemma}
\newtheorem{cor}[theorem]{Corollary}
\newtheorem{corollary}[theorem]{Corollary}
\newtheorem{observation}[theorem]{Observation}
\newtheorem{obs}[theorem]{Observation}
\newtheorem{definition}[theorem]{Definition}
\newtheorem{def*}[section]{Definition}


\def\cA{\mathcal A}
\def\cB{\mathcal B}
\def\cC{\mathcal C}
\def\cD{\mathcal D}
\def\cE{\mathcal E}
\def\cF{\mathcal F}
\def\cG{\mathcal G}
\def\cH{\mathcal H}
\def\cJ{\mathcal J}
\def\cI{\mathcal I}
\def\cK{\mathcal K}
\def\cL{\mathcal L}
\def\cM{\mathbb M}
\def\cN{\mathcal N}
\def\cO{\mathcal O}
\def\cP{\mathcal P}
\def\cQ{\mathcal Q}
\def\cR{\mathcal R}
\def\cS{\mathcal S}
\def\cT{\mathcal T}
\def\cU{\mathcal U}
\def\cV{\mathcal V}
\def\cW{\mathcal W}
\def\cX{\mathcal X}
\def\cY{\mathcal Y}
\def\cZ{\mathcal Z}


\def\OO{\mathrm  O}
\def\Sp{\mathrm  {Sp}}
\def\SO{\mathrm  {SO}}



\def\frA{\mathfrak A}
\def\frB{\mathfrak B}
\def\frC{\mathfrak C}
\def\frD{\mathfrak D}
\def\frE{\mathfrak E}
\def\frF{\mathfrak F}
\def\frG{\mathfrak G}
\def\frH{\mathfrak H}
\def\frJ{\mathfrak J}
\def\frK{\mathfrak K}
\def\frL{\mathfrak L}
\def\frM{\mathfrak M}
\def\frN{\mathfrak N}
\def\frO{\mathfrak O}
\def\frP{\mathfrak P}
\def\frQ{\mathfrak Q}
\def\frR{\mathfrak R}
\def\frS{\mathfrak S}
\def\frT{\mathfrak T}
\def\frU{\mathfrak U}
\def\frV{\mathfrak V}
\def\frW{\mathfrak W}
\def\frX{\mathfrak X}
\def\frY{\mathfrak Y}
\def\frZ{\mathfrak Z}

\def\fra{\mathfrak a}
\def\frb{\mathfrak b}
\def\frc{\mathfrak c}
\def\frd{\mathfrak d}
\def\fre{\mathfrak e}
\def\frf{\mathfrak f}
\def\frg{\mathfrak g}
\def\frh{\mathfrak h}
\def\fri{\mathfrak i}
\def\frj{\mathfrak j}
\def\frk{\mathfrak k}
\def\frl{\mathfrak l}
\def\frm{\mathfrak m}
\def\frn{\mathfrak n}
\def\fro{\mathfrak o}
\def\frp{\mathfrak p}
\def\frq{\mathfrak q}
\def\frr{\mathfrak r}
\def\frs{\mathfrak s}
\def\frt{\mathfrak t}
\def\fru{\mathfrak u}
\def\frv{\mathfrak v}
\def\frw{\mathfrak w}
\def\frx{\mathfrak x}
\def\fry{\mathfrak y}
\def\frz{\mathfrak z}

\def\frsp{\mathfrak{sp}}

\def\frheis{\mathfrak{heis}}


\def\bfa{\mathbf a}
\def\bfb{\mathbf b}
\def\bfc{\mathbf c}
\def\bfd{\mathbf d}
\def\bfe{\mathbf e}
\def\bff{\mathbf f}
\def\bfg{\mathbf g}
\def\bfh{\mathbf h}
\def\bfi{\mathbf i}
\def\bfj{\mathbf j}
\def\bfk{\mathbf k}
\def\bfl{\mathbf l}
\def\bfm{\mathbf m}
\def\bfn{\mathbf n}
\def\bfo{\mathbf o}
\def\bfp{\mathbf p}
\def\bfq{\mathbf q}
\def\bfr{\mathbf r}
\def\bfs{\mathbf s}
\def\bft{\mathbf t}
\def\bfu{\mathbf u}
\def\bfv{\mathbf v}
\def\bfw{\mathbf w}
\def\bfx{\mathbf x}
\def\bfy{\mathbf y}
\def\bfz{\mathbf z}

\def\bfA{\mathbf A}
\def\bfB{\mathbf B}
\def\bfC{\mathbf C}
\def\bfD{\mathbf D}
\def\bfE{\mathbf E}
\def\bfF{\mathbf F}
\def\bfG{\mathbf G}
\def\bfH{\mathbf H}
\def\bfI{\mathbf I}
\def\bfJ{\mathbf J}
\def\bfK{\mathbf K}
\def\bfL{\mathbf L}
\def\bfM{\mathbf M}
\def\bfN{\mathbf N}
\def\bfO{\mathbf O}
\def\bfP{\mathbf P}
\def\bfQ{\mathbf Q}
\def\bfR{\mathbf R}
\def\bfS{\mathbf S}
\def\bfT{\mathbf T}
\def\bfU{\mathbf U}
\def\bfV{\mathbf V}
\def\bfW{\mathbf W}
\def\bfX{\mathbf X}
\def\bfY{\mathbf Y}
\def\bfZ{\mathbf Z}

\def\bfw{\mathbf w}

\def\R {{\mathbb R }}
 \def\C {{\mathbb C }}
  \def\Z{{\mathbb Z}}
  \def\H{{\mathbb H}}
\def\K{{\mathbb K}}
\def\N{{\mathbb N}}
\def\Q{{\mathbb Q}}
\def\A{{\mathbb A}}

\def\T{\mathbb T}
\def\L{\mathbb L}

\def\bbA{\mathbb A}
\def\bbB{\mathbb B}
\def\bbD{\mathbb D}
\def\bbE{\mathbb E}
\def\bbF{\mathbb F}
\def\bbG{\mathbb G}
\def\bbI{\mathbb I}
\def\bbJ{\mathbb J}
\def\bbL{\mathbb L}
\def\bbM{\mathbb M}
\def\bbN{\mathbb N}
\def\bbO{\mathbb O}
\def\bbP{\mathbb P}
\def\bbQ{\mathbb Q}
\def\bbS{\mathbb S}
\def\bbT{\mathbb T}
\def\bbU{\mathbb U}
\def\bbV{\mathbb V}
\def\bbW{\mathbb W}
\def\bbX{\mathbb X}
\def\bbY{\mathbb Y}


\def\sfA{\mathsf{A}}
\def\sfBC{\mathsf{BC}}
\def\sfD{\mathsf{D}}
\def\sfE{\mathsf{E}}
\def\sfF{\mathsf{F}}
\def\sfG{\mathsf{G}}
\def\sfH{\mathsf{H}}
\def\sfR{\mathsf{R}}

\def\G{\Gamma}
\def\bS{\bbS}
\def\bL{\L}

%
 \def\ov{\overline}
\def\wt{\widetilde}
\def\wh{\widehat}

%

\def\obmanka{$\mathop{}\quad$}

\def\le{\leqslant}
\def\ge{\geqslant}
\def\la{\langle}
\def\ra{\rangle}

\def\Section{\S}
\def\Sections{\S\S}

\def\SS{\smallskip}
\def\MS{\medskip}
\def\BS{\bigskip}

\def\phi{\varphi}
\def\epsilon{\varepsilon}
\def\kappa{\varkappa}

\renewcommand{\Re}{\mathop{\rm Re}\nolimits}

\renewcommand{\Im}{\mathop{\rm Im}\nolimits}

\def\codim{\mathop{\mathrm{codim}}\nolimits}

\def\Deck{\mathop{\mathrm{Deck}}\nolimits}


\newfont{\weird}{cmff10}



\SS


\SS


\SS


\SS

%


\title{
Rolling of Coxeter polyhedra along  mirrors}
\author[Alekseevski, Michor, Neretin]{Dmitri V. Alekseevski,
Peter W. Michor, Yurii A. Neretin}
\thanks{P.W.M. was supported
by `Fonds zur F\"orderung der wissenschaftlichen Forschung, Projekt
P~14195~MAT';
\newline
 Yu.A.N was supported by Austrian
``Fonds zur F\"orderung der wissenschaftlichen Forschung'',
 project 19064,
and also by
the Russian Agency on Nuclear Energy,
 the Dutch fund NWO, grant 047.017.015,
and the Japan--Russian grant JSPS--RFBR 07-01-91209.}


\begin{abstract}
The topic of the paper
are  developments of  $n$-dimensional
 Coxeter polyhedra. We show that the
 surface of such polyhedron admits a canonical cutting
 such that each piece can be covered by a Coxeter
 $(n-1)$-dimensional domain.
\end{abstract}

\address{D. V.  Alekseevskij:
School  of Mathematics, Edinburgh University, Edinburgh,  
EH9  3JZ, United Kingdom}
\email{d.aleksee@ed.ac.uk}

\address{P.\  W.\  Michor: Fakult\"at f\"ur Mathematik, Universit\"at Wien,
Nordbergstrasse 15, A-1090 Wien, Austria; {\it and}: Erwin
Schr\"odinger Institute of Mathematical Physics, Boltzmanngasse 9,
A-1090 Wien, Austria.}
 \email{Peter.Michor@esi.ac.at}

\address{Yu. A. Neretin: Fakult\"at f\"ur Mathematik, Universit\"at Wien,
Nordbergstrasse~15, A-1090 Wien, Austria {\it and} Group of
Math.Physics, ITEP, B.Cheremushkinskaya, 25, Moscow, 117259, Russia
,{\it and}  MechMath Department, Moscow State University,
Vorob'yovy Gory, Moscow, Russia}
\email{neretin@mccme.ru}

\maketitle


\section{Introduction. Coxeter groups}
\label{s:preliminaries-1}


\subsection{Coxeter groups in spaces of constant curvature.}
\label{ss:constant-curvature} Consider a Riemannian
space $\bbM^n$ of  constant curvature, i.e., a Euclidean
space $\R^n$, a sphere $\bbS^{n-1}$, or a Lobachevsky
space $\L^n$ (on  geometry of such spaces, see \cite{AVS}).

Let $C\subset\bbM^n$ be
 an intersection of a finite or locally finite collection
of half-spaces%
\footnote{A natural example with an infinite collection
of half-spaces is given on Figure \ref{fig:last}.}.

Consider  reflections of $C$
with respect to all $(n-1)$-dimensional faces.
Next, consider ``new polyhedra'' and their reflections
with respect to their faces. Etc. The domain $C$
 is said to be a {\it Coxeter domain}
if we get a tiling of the whole space in this way.
The group of isometries generated by all such reflections
is said to be a {\it reflection group} or a {\it Coxeter group}
(in a narrow sense, see below).
We say that a Coxeter group is {\it cocompact} if
the initial domain $C$ is compact.
In this case, we  say that $C$ is a
{\it Coxeter polyhedron}.

Evidently, if $C$ is a Coxeter domain,
then the dihedral angles between two neighboring faces
of $C$ are
of the form $\frac{\pi}{m}$, where $m\ge 2$ is an integer.
In particular, they are {\it acute}, i.e., $\le 90^\circ$.

Denote the faces of the polyhedron $C$ by $F_1$, \dots, $F_p$,
denote by $s_1$, \dots, $s_p$ the corresponding reflections.
Denote by $\pi/m_{ij}$ the angles between adjacent faces.
Evidently,
 \begin{equation}
s_j^2=1,\qquad (s_i s_j)^{m_{ij}}=1
.
\label{eq:coxeter-def}
\end{equation}

\subsection{More terminology.}
 Consider a Coxeter
tiling of $\bbM^n$. Below
a ''{\it chamber}'' is any ($n$-dimensional)
 polyhedron of the tiling.
 A ``{\it face}'' or ``{\it facet}'' is
an $(n-1)$-dimensional face of some chamber;
 a {\it hyperedge}
is an $(n-2)$-dimensional edge; a {\it stratum}
is an arbitrary stratum of $\codim\ge 1$
 of some chamber; a {\it vertex} is a vertex.

Also  ``{\it mirrors}'' are  hyperplanes
of  reflections. They divide
the space $\bbM^n$ into chambers.
The group $G$ acts on the set of chambers
simply transitively.
We denote the reflection with respect to a mirror
$Y$ by $s_Y$.

 Each facet is contained in a unique mirror.


\subsection{General Coxeter groups.}\label{Coxeterscheme} Take a symmetric $p\times p$
 matrix $M$ with positive integer elements,
 set $m_{jj}=1$;
 we admit $m_{ij}=\infty$. An
 {\it abstract Coxeter group} is a group
 with generators $s_1$, \dots, $s_n$ and relations
 (\ref{eq:coxeter-def}).

For such a group we draw a graph (we use
the term ``{\it Coxeter scheme}'')
 in the following way. Vertices of the graph
correspond to generators. We connect $i$ and $j$-th vertices
 by $(m_{ij}-2)$ edges.
In fact, we draw a multiple edge if $k\le 6$, otherwise
we write a number $k$ on the edge.

This rule also assign a graph to each Coxeter polyhedron.



\subsection{Spherical Coxeter groups.} By definition,
a spherical Coxeter group, say $\Gamma$,
 acts by orthogonal transformations
of the Euclidean space $\R^{n+1}$. A group
$\Gamma$ is said to be
 {\it reducible} if there exists a proper
  $\Gamma$-invariant subspace in $\R^{n+1}$.
Evidently, the orthogonal complement to a $\Gamma$-invariant
subspace is $\Gamma$-invariant.

The classification of irreducible Coxeter groups
is well known%
\footnote{Actually, these objects were known to Ludwig Schl\"afli
and  Wilhelm Killing in XIX century.
 In 1924, Hermann Weyl identified
these groups as reflection groups,
in 1934 Harold Coxeter  gave a formal classification
and also classified Euclidean groups.}%
, see Bourbaki \cite{Bou}.
The list consists of Weyl groups of semisimple
Lie algebras ($=$ Killing's list of root systems)
$+$  dihedral groups $+$  groups
of symmetries of the icosahedron and  4-dimensional
hypericosahedron (the table is given Section 3).

This also gives a classification of reducible groups.


\subsection{Coxeter equipments.} Next, consider an arbitrary
Coxeter polyhedron in $\R^n$, $\bbS^n$, or $\L^n$.
Consider a stratum $H$ of codimension $k$,
it is an intersection of $k$ faces,
$H=F_{i_1}\cap\dots\cap F_{i_k}$.
The reflections with respect to the faces
$F_{i_1}$, \dots, $F_{i_k}$
generate a Coxeter group, denote it by
$\Gamma(H)=\Gamma(F_{i_1},\dots, F_{i_k})$.

This group is a spherical Coxeter group. Namely,
for $x\in H$ consider the orthocomplement in the tangent space at $x$
to the stratum $H$ and the sphere in this orthocomplement. Then
$\Gamma(H)$ is a reflection group
of this Euclidean sphere.

If $H\subset H'$, then we have the tautological
embedding
$$
\iota_{H',H}:\Gamma(H')\to\Gamma(H)
.
$$
If $H\subset H'\subset H''$,
then
$$
\iota_{H'',H}=\iota_{H',H}\iota_{H'',H'}
.
$$
Such a collection of groups and homomorphisms
is said to be a {\it Coxeter equipment.}


\subsection{Cocompact Euclidean Coxeter groups.}
Here classification is also  simple and
 well known, see Bourbaki \cite{Bou}.
Any such group $\Gamma$
 contains a normal subgroup $\Z^n$ acting by translations
 and $\Gamma/\Z^n$ is a spherical Coxeter group.


\subsection{Coxeter groups in Lobachevsky spaces.}
We report from Vinberg \cite{Vin2}, Vinberg, Shvartsman, \cite{VS}.
The situation  differs drastically.

\SS

a) Coxeter polygons on Lobachevsky plane are arbitrary $k$-gons
with angles of the form $\pi/m_j$. The sum of exterior angles
must satisfy $\sum(\pi-\pi/m_j)>2\pi$. If $k>5$ this condition
holds automatically. For $k=4$ this excludes  rectangles,
also few triangles are forbidden (in fact, spherical and Euclidean
triangles). A Coxeter $k$-gon with prescribed angles
 depends on $(k-3)$ parameters.

\SS

b) In  dimensions $n>2$ Coxeter polyhedra are rigid.
There are many Coxeter groups in spaces
of small dimensions ($n=3$, 4, 5),
but for $n\ge 30$ there is no Coxeter
group with compact fundamental polyhedron at all.
For $n>996$ there is no Coxeter group of  finite covolume
(Prokhorov, Khovanskii, 1986, see \cite{Khov});
the maximal dimensions of known examples are:
8 for compact polyhedra (Bugaenko),
and 21 for a polyhedron of finite volume
(Borcherds).
For $n=3$ there is  a nice Andreev's description
\cite{And} of all Coxeter
polyhedra, it is given in the following two subsections.

\subsection{Acute angle polyhedra in $\L^3$.}
\label{ss:acute}
First, we recall the  famous (and highly non-trivial)
 Steinitz Theorem (see, e.g., \cite{Lus})
about possible combinatorial structure of convex
polyhedra in $\R^3$.

Since the boundary of  a polyhedron is a topological sphere $S^2$,
edges form a connected graph on the sphere, it divides the sphere
into polygonal domain (we use the term 'face'
for such a domain). There are the following
evident properties of the graph:

\SS

--- each edge is contained in 2 faces;

\SS

--- each face has $\ge 3$ vertices;

\SS

--- the intersection of any pair of faces can by the empty set,
a vertex, or an edge.

\begin{theorem*} {\bf (Ernst Steinitz)}
 Any  graph on the sphere $S^2$
 satisfying the above conditions can be realized
 as a graph of edges of a convex polyhedron.
 \end{theorem*}

Our next question is the existence
of a convex polyhedron in $\L^3$ of a given
combinatorial structure where each dihedral (i.e., between two adjacent faces) angle
is a given {\it acute} angle
('acute' or also 'non-obtuse' means $\le \pi/2$)
There are the following
a priori properties of such polyhedra:

\SS

1) All spatial angles are simplicial, i.e.,
each vertex of the graph is contained in 3 edges.
The angles $\phi_1$, $\phi_2$, $\phi_3$
in a given vertex satisfy
\begin{equation}
\phi_1+\phi_2+\phi_3>2\pi
.\label{eq:sum-angles-sphere}
\end{equation}

%
%

\SS

2) At each vertex, 
the set of all dihedral angles determines all other angles in each face at this vertex
(by the spherical cosine theorem).
A face must be a Lo\-ba\-chev\-sky polygon,
i.e., the sum of its exterior angles must be $\ge 2\pi$.
Since all dihedral angles are acute, angles in each face are also acute.
Therefore our conditions forbid  only rectangles
and some triangles.

\SS

3) The following restriction is non-obvious:
We say that  a {\it $k$-prismatic element}
of a convex polyhedron $C$ is a 
sequence
$$
F_{1},\quad F_{2},\quad \dots,\quad F_{k},\quad
F_{k+1}:=F_1
$$
of  faces such that $F_k$ and $F_{k+1}$
have a common edge, and all triple
intersections $F_i\cap F_j\cap F_k$ are empty.

\begin{lemma*}
{\bf (Andreev)}
For any prismatic element in an acute angle polyhedron,
the sum of exterior
dihedral angles is $>2\pi$.
\end{lemma*}

\begin{figure}
$\epsfbox{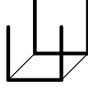}$
\caption{The following configuration with
dihedral  angles
$=\pi/2$ on {\it thick} edges is forbidden
in the Andreev Theorem. In this case,
we would get a quadrangle with right angles, but such quadrangles
do not exist in Lobachevsky space.}
\label{fig:add}
\end{figure}

\begin{theorem*}
{\bf (Andreev)} Consider a Steinitz-admissible
3-valent spherical graph
 with $>4$ vertices%
 \footnote{Simplices are exceptions. However, their
 examination is simple, Lanner, 1950, see, e.g., \cite{VS}.}.
 Prescribe a dihedral acute angle to each edge
 in such a way that:

\SS

---  the inequality
 (\ref{eq:sum-angles-sphere}) in each vertex is satisfied;

 \SS

--- all 3- and 4-prismatic elements satisfy  the previous
 lemma;

\SS

 --- we  forbid the configuration given on
 Figure \ref{fig:add}.

\SS

Under these assumptions, there exists a unique
convex polyhedron $\subset \L^3$
of the given combinatorial structure and
with the given acute angles.
\end{theorem*}

The uniqueness is a rigidity theorem of Cauchy type
(see \cite{Ale},\cite{Lus}).
The existence is a deep unusual fact; it
is a special case of a theorem
of Aleksandrov type \cite{Ale}
obtained by Rivin, see \cite{HR}, \cite{Hod}.

 For some applications of the Andreev and Rivin
 Theorems to elementary  geometry,
 see Thurston \cite{Thu}, Rivin \cite{Riv2}.

\begin{figure}
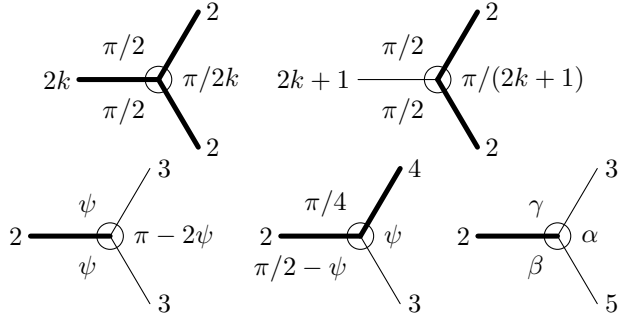

\epsfbox{curvature.1}
\quad
\epsfbox{curvature.2}
\\
\epsfbox{curvature.3}
\quad
\epsfbox{curvature.4}
\quad
\epsfbox{curvature.5}
\caption{We draw all possible types
of vertices of an Andreev polyhedron.
We present the labels $m_j$ on the edges and flat angles
in faces. Here $\psi=\arctan \sqrt 2$ and
$\alpha$, $\beta$, $\gamma$ are explicit angles
with $\alpha+\beta+\gamma=\pi/2$.
Evaluations of all these angles are given in figures in Section \ref{s:Rolling}.
\newline
\obmanka
We draw a {\it thick} line iff the label is even.}
\label{fig:andreev-vertices}
\end{figure}

\subsection{Andreev polyhedra.}
Andreev's Theorem provides us a description of all
Coxeter polyhedra in $\L^3$. Now all angles  have the form
$\pi/m_{ij}$ with integer $m_{ij}>1$.
We  simply write the  labels $m_{ij}$
on the corresponding edges.

Below the term ``{\it Andreev polyhedron}'' will mean
a compact Coxeter polyhedron in $\L^3$.

\SS

All possible pictures at vertices of  Andreev
polyhedra are given in Figure \ref{fig:andreev-vertices}.

\subsection{Results of the paper.}
Consider a convex polyhedron $C$
in a space $\bbM^n$ of constant curvature.
Following Alexandrov \cite{Ale}, we regard
the boundary $\Xi=\partial C$ of $C$ as an $(n-1)$-dimensional
manifold of constant curvature with singularities.
In the case $n=3$, we get a two-dimensional surface
with conic singularities of negative curvature
(see e.g. Figure \ref{fig:andreev-vertices}, in all the cases
the sum of angles at a singularity is $<2\pi$).

Now, cut $\Xi$ along hyperedges with {\it even} labels
(i.e., hyperedges with dihedral angles $\pi/2k$).
Let $\Omega_1$, $\Omega_2$, \dots be the connected pieces
of the cut 
surface.

\begin{theorem}
The universal covering $\Omega_j^\sim$
 of each $\Omega_j$ is a Coxeter domain
in $\bbM^{n-1}$.
\end{theorem}

{\sc Proof for Andreev polyhedra.}
We simply look to Figure  \ref{fig:andreev-vertices}.
In all the cases, angles between {\it thick}
edges are Coxeter.
\hfill $\square$

\SS

We also describe tilings of mirrors, groups of transformations
  of  mirrors induced by the initial Coxeter
group (Theorem \ref{semidirect})
and the  Coxeter equipments of $\Omega_j^\sim$
(Theorem \ref{th:reduction}).

\SS

The addendum to the paper contains two examples of 'calculation'
of developments,  for an Andreev prism $\subset\L^3$
and for a Coxeter simplex $\subset\L^4$.
The proof of the Andreev Theorem is nonconstructive.
In various explicit cases, our
argumentation  allows to construct
an Andreev polyhedron from the a priori information about
its development. Our example illustrates this phenomenon.

On the other hand, there arises a natural problem of elementary
geometry:

\SS

--- {\it Which Andreev polyhedra are partial developments
of 4-dimensional
Coxeter polyhedra? Is it possible to describe all 3-dimensional
polyhedra that are faces of 4-dimensional Coxeter polyhedra?}

\SS

Our main argument (Rolling Lemma \ref{rolling-lemma})
is very simple, it is valid in a wider generality,
we briefly discuss such possibilities in the next two subsections.

\subsection{Polyhedral complexes and projective Coxeter polyhedra.}
\label{ss:projective-coxeter}

\begin{theorem*}{\bf(Tits)}
 Any Coxeter group can be realized as a group of transformations of an open convex 
 subset of a real projective space  $\R \bbP^n$ which is
generated by a collection of reflections
$s_1$, \dots, $s_p$ with respect to  hyperplanes\footnote{a reflection is determined by a 
fixed hyperplane and a reflected tranversal line}
 intersecting the subset. 
The closure of a chamber is a convex polyhedron.
\end{theorem*}

See also Vinberg \cite{Vin}.


\subsection{A more general view.}
\label{ss:wider}
Nikolas Bourbaki\footnote{Apparently,
he used the work by Jacque Tits \cite{Tit};
the latter text is inaccessible  for the authors.}
 proposed a way to build topological spaces
 from  Coxeter groups.
M.~Davis used this approach
in numerous papers (see e.g. \cite{Dav}, \cite{Dav2})
and the book \cite{Dav-book}; in particular he
constructed nice  examples/counterexamples in topology.

Also it is possible to consider arbitrary Riemannian
manifolds equipped with a discrete isometric action
of a Coxeter group such that the set of fixed
points of each generator is a (totally geodesic) hypersurface, and such that
the generators act as reflections
with respect to these submanifolds. In this context,
a chamber itself can be a topologically non-trivial
object, see \cite{Dav-book},
\cite{reflec}.

\section{Rolling of chamber}
\label{s:Rolling}

In this section,
  $\bbM^n$ is a  space
  $\L^n$, $\bbS^n$, $\R^n$ of constant curvature
  equipped with a Coxeter group
$\Gamma$ or, more generally, any space described in
Subsection
\ref{ss:projective-coxeter}.

\begin{figure}
\epsfbox{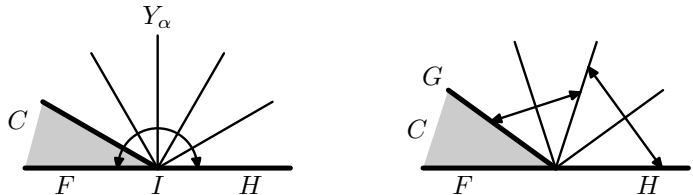}
\caption{Even and odd labels. Proof of Rolling Lemma.}\label{current.1}
\end{figure}

Fix a mirror $\bbX^{n-1}$ in  $\bbM^n$.
Consider intersections of $\bbX^{n-1}$ with
other mirrors $Y_\alpha$.
The set $\bbX^{n-1}\setminus \bigcup Y_\alpha$ is a disjoint
union of open facets.
Thus, we get a tiling of $\bbX^{n-1}$ by facets.

Our aim is to
describe this tiling in the terms of the geometry of a chamber.

\subsection{Rolling lemma}
\label{ss:rolling}

\begin{lem}
\label{rolling-lemma}
Let $I\subset \bbX^{n-1}$ be an
$(n-2)$-dimensional hyper-edge of our tiling.
Let $F$, $H\subset \bbX^{n-1}$ be the facets adjacent  to $I$.

\SS

a)      If the label $m_I$ of $I$ is even, then
$I$ is contained in a certain mirror $Y_\alpha$ orthogonal
to $X$. In particular $s_{Y_\alpha} F=H$.

\SS

b)      Let the label $m$ be odd. Let $C$ be a chamber
adjacent to the facet $F$. Let $G$ be another facet of $C$
adjacent to the same hyper-edge $I$. Then $G$ is isometric to
$H$. More precisely, there is $\gamma\in\G$ fixing all the points
of $I$
such that $\gamma G=H$.
\end{lem}

{\sc Proof} is given in figure \ref{current.1}.

\begin{figure}
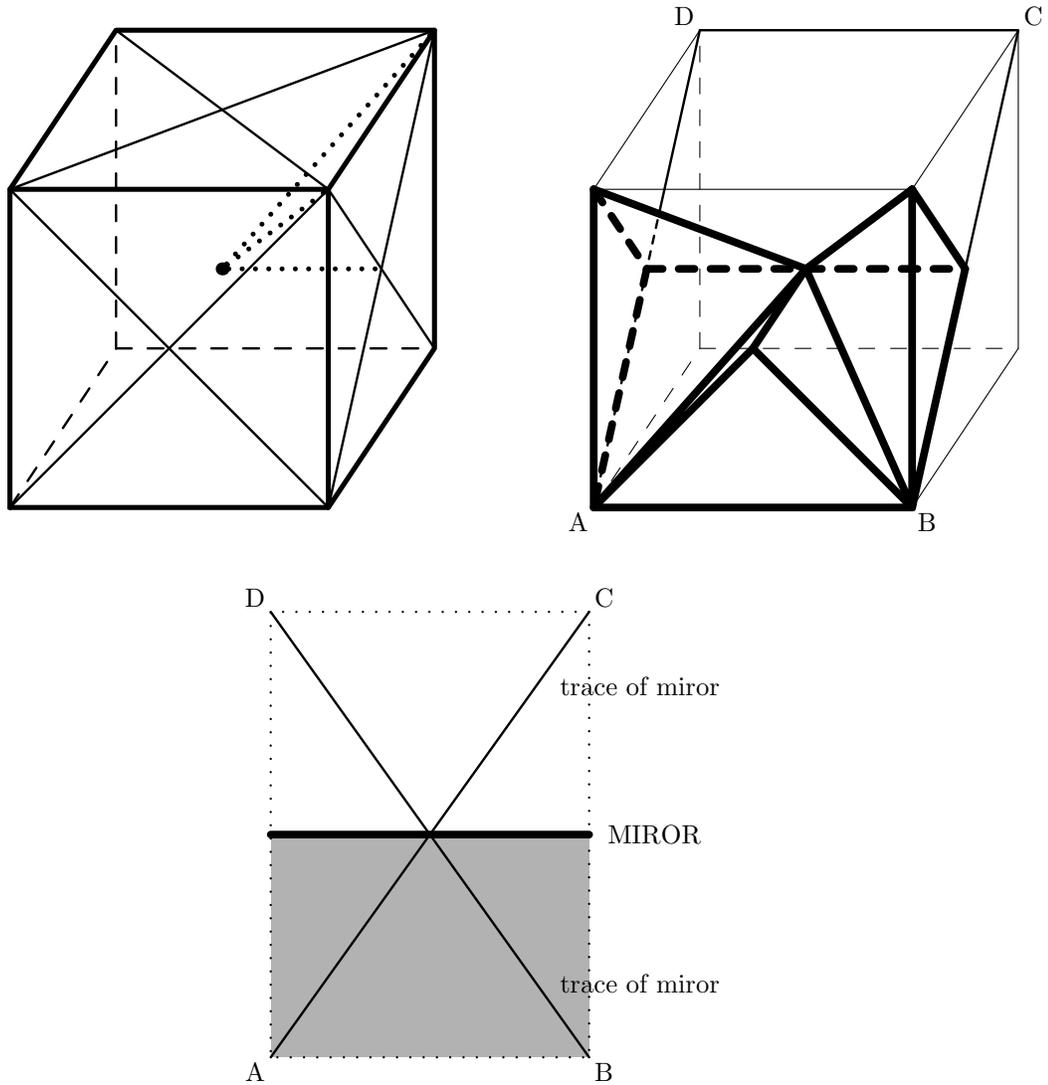

$$\epsfbox{cube.1}$$

$$\epsfbox{cube.2}$$
\caption{Example of rolling:
 the reflection group $\sfA_3$
in $\R^3$. The mirrors are planes passing through
opposite edges of the cube.
There are 24 Weyl chambers, which are simplicial cones
with dihedral angles $\pi/3$, $\pi/3$, $\pi/2$
(we draw them as simplices). Rolling
of a Weyl chamber by the mirror $ABCD$ produces a half-plane.
\newline
\obmanka
We can also  regard $\sfA_3$ as a reflection group
on the 2-dimensional sphere $\bbS^2$.}\label{cube.12}
\end{figure}

\begin{figure}
$$\epsfbox{icos.1}$$
\caption
{ Example of rolling:
the icosahedral group $\sfH_3$. It is generated by reflections
with respect to bisectors of segments connecting
midpoints of opposite edges of the icosahedron.
The bisectors separate $\R^3$ into
120 simplicial cones
with dihedral angles $\pi/2$, $\pi/3$, $\pi/5$.
In the figure the simplicial
cones are cut by the surface of the icosahedron.
\newline
\obmanka
We show an admissible rolling of
a Weyl chamber along a mirror.
The final chamber in the mirror is a quadrant.}
\label{icos}

\end{figure}


%
%

\begin{figure}
$$
\epsfbox{cube.3}
$$
\caption{ Example of rolling:
 the (affine) Euclidean
  reflection group $\wt A_4$ in $\R^3$.
A chamber is the simplex $ABCD$.
Rolling through $AB$ and $CD$ is forbidden.
Deleting these edges from the surface of the simplex, we get
a non-simply  connected surface.
Hence, the process of rolling  is infinite.
The arrow shows the deck transformation induced by the generator
of the fundamental group.}\label{cube.3}


%
%

\bigskip

$$
\epsfbox{cube.4}
$$
\caption{Example of rolling:
 the hypericosahedral group $\sfH_4$
 acting
on the 3-dimensional sphere $\bS^3$. The chamber is
is the spherical simplex drawn
in the figure (we omit all labels $\pi/2$ on edges).
\newline
\obmanka
The angle $=\pi$ on the development at
 $D$ was
evaluated in figure \ref{cube.12}. The right angle
at $C$ was evaluated in figure \ref{icos}.
\newline
\obmanka
The spherical triangle $ABC$ is present on the
circumscribed sphere in the next figure
(in spite of the absence of the sphere itself).
}\label{cube.4}
\end{figure}


%
%

\begin{figure}
$$
\epsfbox{ikos.1}
$$
\caption{ Example of rolling:
 the hypericosahedral
 group $\sfH_4$
acting in $\R^4$. The figure presents the tiling
of a mirror, i.e., of $\R^3$, by simplicial cones.
We draw intersections of simplicial cones with
the boundary of the icosahedron.
Consider 3 types of 'axes' of the icosahedron:
\newline
\obmanka
A)      segments connecting midpoints of opposite
edges;
\newline
\obmanka
B) segments connecting central points of opposite faces;
\newline
\obmanka
C) diagonals connecting opposite vertices.
\newline
\obmanka
Consider  bisectors of all such segments.
Type A bisectors are mirrors. They divide $\R^3$
into 120 simplicial chambers. Six chambers are presented in the
front face of the icosahedron.
\newline
\obmanka
Adding bisectors of type B and C
we obtain a partition of $\R^3$ into 480 simplicial cones.
This is the desired tiling.
\newline
\obmanka
In this figure, we present  subdivisions of two chambers.
A proof of this picture is contained on figure \ref{cube.4}
}\label{ikos.1}
\end{figure}

\begin{figure}
$$\epsfbox{cube.5}
\quad
\epsfbox{last.1}
$$
\caption{Example. A Coxeter simplex
in $\L^3$.  Its development is an infinite 'strip'
$\subset \L^2$
bounded by two infinite polygonal curves, interior angles
between segments of polygonal curves are $\pi/2$ and $\pi$.
\label{fig:last}}
\end{figure}


\subsection{Algorithm generating the tiling}\label{algorithm}
Let $C$ be a chamber adjacent to a facet $F\subset \bbX^{n-1}$.
Consider an hyper-edge $I$ of $C$ lying in $\bbX^{n-1}$.

\SS

{\sc Operation 1.} Let the hyper-edge $I$ be odd. Consider a facet $G\ne F$
of $C$ adjacent to $I$, consider the corresponding
$\gamma$ from Lemma  \ref{rolling-lemma}
and draw $\gamma G$ on $\bbX^{n-1}$.

\SS

{\sc Operation 2.} If the hyper-edge $I$ is even,
then we reflect $F$ in $\bbX^{n-1}$ with respect
to $I$.

\SS

We perform all the possible finite sequences of such operations.
By the Rolling Lemma, we get the whole tiling
of the mirror $\bbX^{n-1}$.

\SS

{\sc Remark.} Let $\bbM^n=\R^3$, $\bS^3$, $\bL^3$
be a usual 3-dimensional
space of constant curvature.
Operation 1 corresponds to rolling
of a polyhedron $C$ along the hyperplane
$\bbX^{n-1}\sim \bbM^2$ over the edge $I$.
\hfill$\square$

\subsection{The group preserving the mirror
 $\bbX^{n-1}$.}
 \label{mirror-group}
For a mirror $\bbX^{n-1}$,
consider the group $\G_*=\G_*(\bbX^{n-1})$ of all the isometries
of $\bbX^{n-1}$ induced by elements of $\G$ preserving $\bbX^{n-1}$.

If $\gamma\in \G$
preserves $\bbX^{n-1}$, then $s_{\bbX^{n-1}}\gamma$
also preserves $\bbX^{n-1}$ and agrees with $\gamma$ on $\bbX^{n-1}$.
Thus each element of $\G_*$ is
induced by two different elements of $\G$.

\begin{observation}
Let $F_1, F_2\subset \bbX^{n-1}$ be equivalent facets.
There is a unique element $\mu\in\G_*(\bbX^{n-1})$ such that
$\mu F_1=F_2$.
\end{observation}

\subsection{Reflections in mirrors and the new chamber.}
\label{mirrors-on-mirrors}
Consider all the mirrors $Z_\alpha\subset\bbM^n$ orthogonal
to our mirror $\bbX^{n-1}$. The corresponding
 reflections $s_{Z_\alpha}$
generate a reflection group on
$\bbX^{n-1}$; denote this group by $\Delta=\Delta(\bbX^{n-1})$.

\begin{observation}
 $\Delta$ is a normal subgroup
in $\G_*(\bbX^{n-1})$.
\end{observation}

Indeed, if $s$ is a reflection, then
$\gamma^{-1} s\gamma$ is a reflection.
\hfill$\square$

\SS

Consider a chamber $C$ of $\bbM^n$ lying on
 $\bbX^{n-1}$
(i.e., having a
facet in $\bbX^{n-1}$)
and consider all possible sequences of admissible
rollings, i.e., we allow Operation 1 of Algorithm \ref{algorithm}
and we forbid Operation 2.
Denote by $B\subset \bbX^{n-1}$
the domain obtained by rolling, tiled by the traces of facets of $C$ making contact with 
$\bbX^{n-1}$ during rolling.

\begin{thm}\label{chamber-on-mirror}
$B$ is a chamber of the reflection group $\Delta(\bbX^{n-1})$.
\end{thm}


{\sc Proof.} We can not roll further
if and only if we meet a ``vertical'' mirror.
\hfill $\square$

\SS

{\sc Examples of rolling.}
Some  examples of rolling corresponding to
the usual spherical Coxeter groups

$\sfA_3$ : $
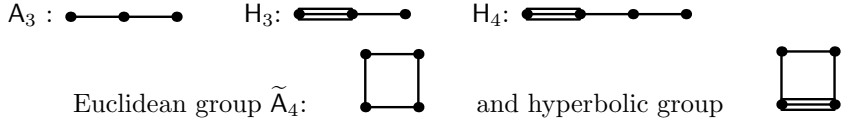
\epsfbox{current.5}$
\qquad
 $\sfH_3$: $\epsfbox{current.8}$
\qquad
$\sfH_4$: $\epsfbox{current.7}$
\qquad
$$
\text{
Euclidean group
$\wt \sfA_4$:}\qquad \epsfbox{current.6}\qquad
\text{and hyperbolic group}\qquad
\epsfbox{current.10}
$$

\SS

\noindent
 are given in 
figure \ref{cube.12}-\ref{fig:last}. In these figures,
we also evaluate the new chamber $B$.\hfill $\square$

\begin{lem}
\label{lemma-orthogonal-mirrors}
Each $(n-3)$-dimensional stratum
of our tiling of $\bbX^{n-1}$ is contained in a mirror
of the group $\Delta(\bbX^{n-1})$.
\end{lem}

{\sc Proof.} This stratum is equipped  with
a finite 3-dimensional Coxeter group
(i.e., $\sfA_3$, $\sfBC_3$, $\sfH_3$, $\sfA_1\oplus \sfG_2^{m}$,
$\sfA_1\oplus \sfA_1\oplus \sfA_1$,
see Table below).
For each mirror of such a group there exists  an orthogonal
mirror.
\hfill$\square$

\subsection{Rolling scheme.}\label{rolling-scheme}
Denote by $\Xi( C)$ the surface of the
initial  chamber $C$, let
$\Xi'(C )$ be the surface with all even edges deleted.

\begin{lem}
 $\Xi'$ does not contain  $(n-3)$-dimensional strata
of $C$.
\label{no-vertices}
\end{lem}

This is rephrasing of  Lemma
 \ref{lemma-orthogonal-mirrors}. \hfill $\square$

\SS

Consider the graph, whose vertices are the
facets of $\Xi'$; vertices are connected by an edge
if the corresponding facets
are neighbors in $\Xi'$.
We call this graph the {\it Rolling scheme}.
In fact, the Rolling scheme is the Coxeter scheme \ref{Coxeterscheme}
with removed even (and infinite) edges.

\begin{prop}
\label{p:homotopy}
The surface
  $\Xi'$ is homotopically
equivalent to the Rolling scheme.
\end{prop}

\subsection{Proof of Proposition \ref{p:homotopy}.}
Let $U$ be a convex polyhedron in $\R^n$,
denote by $\Xi$ its surface.
Choose a point $A_j$ in interior of each
$(n-1)$-dimensional face.
Choose
a point $B_k$ in interior of each $(l-2)$-dimensional
boundary stratum (hyperedge) of $U$.

Draw the segment $[A_j,B_k]$ iff the face contains the hyperedge.
Thus we get a  graph $T$
 on the surface of the polyhedron $C$
 whose vertices are enumerated by faces of $U$ and edges
 are enumerated by hyperedges of $U$.
Denote by $\Xi^\bigtriangledown$
the surface of the polyhedron $S$ without
boundary strata of dimension $(n-3)$.

\begin{lem}
\label{l:retraction}
 The graph $T$ is a deformation
retract of  $\Xi^\bigtriangledown$.
Moreover, it is possible to choose a homotopy
that preserves all faces
and all hyperedges.
 \end{lem}

{\sc Proof.} See figure \ref{fig:retraction}.
\hfill $\square$

\SS

\begin{figure}
a) \epsfbox{homotopy.1}
b)\quad  \epsfbox{homotopy.2}
\caption{Proof of Lemma \ref{l:retraction}.
\newline
\obmanka
a) $n=3$. Graph on a surface of a 3-dimensional
polytop and a retraction. Recall that we
 have removed vertices.
 \newline
 \obmanka
b) $n=4$. A piece of a 3-face of 4-dimensional polyhedron.
Recall that 1-dimensional edges are removed.
Inside a simplex $P_1 P_2 A_j B_k$ the retraction is the projection
to $A_jB_k$ with center on the segment $P_1P_2$. Note
that all segments connecting $A_jB_k$ and $P_1P_2$ are pairwise
non-intersecting.}
 \label{fig:retraction}
\end{figure}

 Proposition \ref{p:homotopy}
follows from  Lemma \ref{l:retraction}. \hfill $\square$



%
%
%
%
%

\subsection{Action of the
fundamental group on mirror.}\label{pi1-mirror}
Let $F$ be a facet in $\bbX^{n-1}$, let $C$ be a chamber
of $\bbM^n$ lying on $F$, and let $B\supset F$ be the chamber
of the reflection
group
 $\Delta(\bbX^{n-1})$ obtained by rolling $C$, as described in
subsection \ref{mirrors-on-mirrors}.

Let $\Omega$ be a connected component
of $\Xi'$ containing the facet $F$.

Let $F_1$, \dots, $F_r$ be facets $\subset \Omega$.
We can think that each facet has its own color;
thus the mirror $\bbX^{n-1}$ is painted in $r$ colors.
Moreover, for each facet $H\in\bbX^{n-1}$ there is a
canonical bijection ('parametrization') from the corresponding
 $F_i\subset \Omega$ to $H$.
We say that a bijection $\bbX^{n-1}\to\bbX^{n-1}$
(or $B\to B$) is an isomorphism if it
preserves the coloring and commutes with
the parameterizations.

\begin{prop}\label{deck-transformations}
a)      The chamber $B\subset \bbX^{n-1}$ is
 the universal covering of $\Omega$.

\SS

b)      Any deck transformations of $B$ is an isomorphism
$B\to B$
and admits a unique extension to an isomorphism of the mirror
$\bbX^{n-1}$.

\SS

c)      Each isomorphism
 $\mu\in\G_*(\bbX^{n-1})$ preserving $B$ is induced
by a deck transformation.
\end{prop}

{\sc Proof.} a) Denote by $\Omega^\sim$ the universal covering
of $\Omega$. The chamber $B$ was constructed as the image
of $\Omega^\sim$. Moreover,  the map $\Omega^\sim\to C$
is locally bijective. On the other hand, a chamber
on a simply connected manifold is simply connected
see (see \cite{reflec}, 2.14);
therefore $B\simeq\Omega^\sim$.

\SS

b) A deck transformation $B\to B$
 is an isometry by the rolling rules.
 Let a deck transformation send a facet $F$ to $F'$.
Then the facets $F$, $F'$ are $\Gamma$-equivalent,
and the corresponding
map in $\Gamma$ is an isometry of $\bbX^{n-1}$.

\SS

c) Let $F\subset \bbX^{n-1}$ be a facet.
We take the deck transformation
sending $F$ to $F'$.

\subsection{Description of $\G_*(\bbX^{n-1})$.}
\label{description}

\begin{thm}\label{semidirect}
The group $\G_*(\bbX^{n-1})$ is a semidirect product
$\Deck(B)\ltimes \Delta(\bbX^{n-1})$.
\end{thm}

{\sc Proof.} Indeed, the group $\Delta(\bbX^{n-1})$
acts simply transitively
on the set of chambers in $\bbX^{n-1}$;
the group $\Deck(B)$ acts simply
transitively on the set of facets of a given type
in the chamber $B$.
\hfill $\square$

\begin{figure}
\epsfbox{erunda.1}
\caption{A graph of vicinity of $(n-1)$-dimensional facets
in the new $(n-1)$-dimensional chamber $B$ is a tree.}
\end{figure}

\begin{figure}
$$
\epsfbox{current.9}
$$
\caption{Subdivision of the cone normal to a stratum.}
\label{current.9}
\end{figure}


\begin{figure}

$$\text{Table. Reduction of spherical Coxeter schemes}$$

\begin{align*}
&\sfA_n:&\quad&\epsfbox{dynkin.1}& \qquad &\mapsto\qquad \sfA_{n-2}\oplus\R
\\
&\sfBC_n:&\quad&\epsfbox{dynkin.2}& \qquad &\mapsto\qquad
      \sfA_1\oplus \sfD_{n-2}\quad\text{or}\quad \sfBC_{n-1}
\\
&\sfD_n:&\quad&\epsfbox{dynkin.3}& \qquad &\mapsto\qquad \sfA_1 \oplus \sfD_{n-2}
\\
&\sfE_6:&\quad&\epsfbox{dynkin.4}& \qquad &\mapsto\qquad \sfA_5
\\
&\sfE_7:& \qquad&\epsfbox{dynkin.5}& \qquad &\mapsto\qquad \sfD_6
\\
&\sfE_8:& \qquad&\epsfbox{dynkin.6}& \qquad &\mapsto\qquad \sfE_7
\\
&\sfF_4:& \qquad&\epsfbox{dynkin.7}& \qquad &\mapsto\qquad \sfBC_3
       \quad\text{or}\quad \sfBC_3
\\
&\sfG_2^{(m)}:& \qquad &\epsfbox{dynkin.8}&\qquad&\mapsto
               \begin{cases} \text{ $\sfA_1$ or $\sfA_1$ if $m$ is even}
               \\
              \text{$\sfR$, if $m$ is odd}
               \end{cases}
\\
&\sfH_3:& \qquad &\epsfbox{dynkin.10}&\qquad &\mapsto \sfA_1\oplus \sfA_1
\\
&\sfH_4:&\qquad &\epsfbox{dynkin.9}&\qquad& \mapsto \sfH_3
              \end{align*}
\end{figure}

\section{Reduction of equipment}
\label{s:equipment}

We keep the notation of the previous section.
Our aim is to describe  the Coxeter equipment of the new chamber $B$.

\subsection{Combinatorial structure of the tiling
of the chamber.}
Consider a graph $\frF$ whose vertices are enumerated by
$(n-1)$-facets lying in $B$, two vertices are connected by
an edge if they have a common $(n-2)$-dimensional
stratum (a former hyperedge in $\bbM^n$).

\begin{observation}
$\frF$ is a tree.
\end{observation}

{\sc Proof.} Indeed, the universal covering of
a graph is a tree. \hfill $\square$

\SS

If the initial rolling scheme  is a tree,
then we get the same
tree. If the rolling scheme contains a cycle,
then we  get an infinite tree (examples: Figures \ref{cube.3},
\ref{fig:last}, the rolling schemes contain 1 cycle).


\subsection{New equipment.}
 All the strata
of $B$ of
dimension $<(n-2)$ are contained in the boundary of $B$.
These strata of $B$ have their
own equipments
(in the sense of the Coxeter manifold $\bbX^{n-1}$).

For a  boundary stratum $H$ of $B$ and some point $y\in H$, 
denote by $N_H\subset T_y\bbX^{n-1}$ the
normal subspace to $H\subset\bbX^{n-1}$.
The {\it normal cone} $D_H\subset N_H$
is the cone consisting of vectors looking inside $B$.
Some of $(n-2)$-dimensional strata (former hyperedges)
$V_\alpha$ contain $H$ and thus we get the
subdivision of the normal cone $D_H$ by tangent spaces
to $(n-2)$-dimensional strata, see figure \ref{current.9}.

We wish to describe the equipment of $B\subset\bbX^{n-1}$
and the subdivisions of normal cones  $D_H$.

\subsection{Finite Coxeter groups}\label{finite-groups}
Let $\G$ be a finite Coxeter group acting in $\R^n$.
Let $\bbX^{n-1}_j$ be the mirrors, let $v_j$ be the vectors orthogonal
to the corresponding mirrors.
 For a vector $v_k$ denote by $R=R_k$
 the set of all $i$ such that
$v_i$ is orthogonal to $v_k$.

 The reflection group $\Delta(\bbX^{n-1}_k)$
is generated by reflections with respect to mirrors
$\bbX^{n-1}_i$, where $i$ ranges in $R$.

\SS

A.      {\it Let the Coxeter group $\Gamma$ be irreducible.}
 We come to the  list given in the Table.
Some comments:

\SS

1)  $\sfG_2^{(m)}$ denotes
the group of symmetries of a regular plane $m$-gon,
$\sfR$ denotes the one-element group acting in $\R^1$;
 all other notations are standard, see \cite{Bou}.

\SS

2) In some cases,  there are two $\G$-nonequivalent mirrors,
then we write both possible variants.

\SS

The rolling scheme (see \ref{rolling-scheme}) is the
Coxeter scheme without even edges.

\smallskip

{\sc Example.} a) For the Weyl chamber $\sfE_8$,
its complete development
is the Weyl chamber $\sfE_7$.

\SS

b) For the Weyl chamber $\sfBC_n$, one of the facets is
the Weyl chamber $\sfBC_{n-1}$. All the remaining facets
are  connected by the rolling graph; the
development is the Weyl chamber $\sfA_1\oplus \sfD_{n-2}$.

\SS

{\sc Proof of Table} is a case-by-case examination of root systems;
for the groups $\sfH_3$ and $\sfH_4$ the proofs are given Figures 
\ref{icos}, \ref{cube.4}, \ref{ikos.1}
(on the other hand the reader can find a nice coordinate description
of the hypericosahedron in \cite{Bou}.

\SS

B. {\it If the Coxeter group $\G$ be reducible,}
$$\G=\G_1\times \G_2\times\dots $$
then its Weyl chamber is the product of the Weyl chambers
for the corresponding chambers $C=C_1\times C_2\times\dots$.
The Coxeter scheme of $\G$ is the union of the Coxeter schemes of
$\G_j$, hence the rolling graph of $\G$ is
 the union of the rolling graphs for all the $\G_j$.
Now  we reduce one of factors $C_j\mapsto B_j$
 according to the rules
given in the Table, and we get a Weyl chamber
$B_j\times \prod_{i\ne j} C_i$.

\subsection{Reduction of equipment.}\label{reduction}
\label{ss:equipment}

 Let $H$ be an $(n-k)$-dimensional stratum of $C$
($k\ge 3$), let $\Gamma_H(C)$ be the corresponding Coxeter group,
and let $\frN_H(C)$ be its chamber in the normal cone.
Denote by $\Gamma_H(B)$ the corresponding group
of the equipment of $B$ and by $\frN_H(B)$ the corresponding
chamber in the normal cone.

\begin{theorem}%
\label{th:reduction}
The group $\Gamma_H(B)$ is obtained by reduction of the group
$\Gamma_H(B)$ and the subdivision of $\frN_H(B)$
is a partial development of the Weyl chamber $\frN_H(C)$
\end{theorem}

{\sc Proof} is obvious. We consider rollings of $C$ with fixed
hyperedge $H$.
The subdivision of the cone $D_H$ is obtained by rolling
with respect to the hyperedges containing $H$.
\hfill $\square$

\section{Addendum. Elementary geometry of Andreev polyhedra}
\label{s:andreev}


\begin{figure}
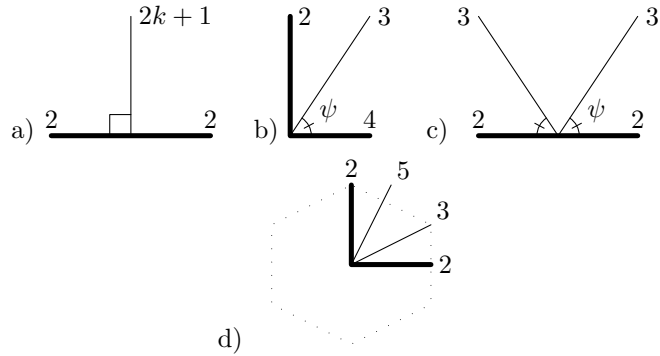

a) \epsfbox{billiard.1}
\quad
b) \epsfbox{billiard.3}
\quad
c) \epsfbox{billiard.2}
$$\text{d)}\quad \epsfbox{billiard.4}$$
\caption{Subdivisions of a Coxeter polygon on the Lobachevsky
plane (we also draw the labels on lines).
There only 4 possible variants of  meetings
between  lines of a subdivision and the boundary.
Here $\tan\psi=\sqrt 2$; for case d) see
figure  \ref{icos}.
\newline
\obmanka
In cases c) and d), the corresponding trihedral angle
of the Andreev polytope is covered by our bended polygon.
}
\label{billiard}
\end{figure}


\begin{figure}
$$\epsfbox{parket.1}$$
\caption{An example of an
Andreev polyhedron in $\bL^3$;
we label the dihedral angles $\pi/3$,
all other dihedral angles are $\pi/2$.
\newline
\obmanka
Its development
 is  a (nonregular) 6-gon, whose angles are $\pi/2$.
The lines $\ell_1$, $\ell_2$ are axes of symmetry.
The polygonal curve $ABDCA$ is a billiard trajectory.
\newline
\obmanka
It is easy to reconstruct the lengths of edges of the prism
 from
the combinatorial structure of the development and the billiard
trajectory.
Indeed, we know the angles of the triangle $AEC$ and of the
``trapezoids``
$ABFE$, and the equiangular
quadrangle $ABDC$.}
\label{parket}
\end{figure}

\subsection{Rolling of Andreev polyhedra and billiard
trajectories in Coxeter polygons.}
Firstly, our construction gives some information
about  developments of  Andreev polyhedra.

Let us roll an Andreev polyhedron
$\subset \bL^3$ along
a mirror $\simeq \L^2$.
In this case, the chamber $B$ of
a mirror is a convex plane Coxeter domain.
By construction, $B$ is subdivided
into several convex polygons
by a certain family of lines.

\begin{prop}
All the possible variants of meetings of
lines of the subdivision  and the boundary of $B$
are presented in figure \ref{billiard}.
\end{prop}

{\sc Proof.} We watch all the possible
variants of reduction of 3-dimensional finite Coxeter
groups to a mirror. The parts a), b), c), d) of figure \ref{billiard}
correspond to $\sfG_2^{2k+1}$,
 $\sfBC_3$,  $\sfA_3=\sfD_3$,  $\sfH_3$,
respectively.
  \hfill $\square$

\SS

\begin{observation}
 The surface of an Andreev polyhedron
is glued from several bended Coxeter polygons;
the rules of bending and the rules of glueing are very rigid.
\end{observation}

{\sc Examples} of rolling of a $3$-dimensional Coxeter polyhedron
are given in figures \ref{fig:last} and \ref{parket}.\hfill $\square$

\begin{figure}
$$\epsfbox{andreev.1}\qquad \epsfbox{andreev.2}$$
\caption{This  prism in $\bL^3$
is a complete development of the Coxeter simplex
$ABCDE$ in $\bL^4$ described in \ref{around}. It carries
all 2-dimensional hyperedges of the initial simplex.
\newline
\obmanka
The development of the prism is the regular 10-gon with
right angles (it also carries all 1-dimensional
strata of the 4-dimensional simplex).
}
\label{andreev}
\end{figure}

\subsection{Example: Rolling along Andreev polyhedra}
\label{around}
Secondly, take a  Coxeter polyhedron in $\bL^4$.
Rolling it
along 3-dimensional Lobachevsky space, we obtain
a Coxeter polyhedron in $\bL^3$ and also some strange
subdivision of this polyhedron.

We present an example.
Consider the simplex $\Sigma$
in $\bL^4$ defined by the Coxeter scheme
\begin{equation}
\epsfbox{dynkin.11}.
\label{dynkin-hyp}
\end{equation}
By $A$,\dots, $E$ we denote the vertices
of the simplex opposite to the corresponding faces. See figure \ref{andreev}.

\SS

{\it Comments to Figure \ref{andreev}.}
The development of $\Sigma$ is a prism drawn in
figure \ref{andreev}.  We write labels for the dihedral angles
 $\ne \pi/2$. Below a ``stratum'' means
a stratum of the tiling; in particular, the
vertical ``edge'' $AB$ consists
 of two 1-dimensional {\it strata} $BC$
and $CA$ and three 0-dimensional ones, $A$, $B$, and $C$.

\SS

1) This is a development. Hence any two strata
 (segments, triangles) having the same
notation are equal (for instance $CD=CD$, $CE=CE$,
$\triangle CBE=\triangle CBE$, etc.).

\smallskip

2) Each stratum (a vertex, a segment) is equipped with
a Coxeter group (this group is visible from its dihedral angles)

\SS

3)   Subdivision of the normal cone $D_H$
to a stratum $H$ (a vertex, a segment)
is determined by the reduction procedure from
subsection \ref{finite-groups}.

\SS

For instance, in the vertex $A$ we have the
subdivision
of the spherical triangle $\sfH_3$ drawn in figure \ref{cube.4},
$$
\epsfbox{cut.1}\qquad\qquad \mapsto \qquad \qquad \sfH_3.
$$

In the normal cone to the
 edge-stratum $DE$ of the prism, we have the icosahedral
  subdivision, see figure \ref{ikos.1},
$$
\epsfbox{cut.2}\qquad\qquad \mapsto \qquad \qquad A_1\oplus A_1.
$$

The normal cone to the segment
$AE$ is drawn in figure \ref{cube.12}; in particular, both
angles
of incidence are $\arctan \sqrt 2$,
$$
\epsfbox{cut.3}\qquad\qquad \mapsto \qquad A_1\oplus \sfR.
$$

The ``front'' face $ABDB$ is orthogonal to the sections
$CDE$ and $ADE$ (since the lines $CD$ and $AD$ of intersection
are equipped with the group $\sfA_1\oplus \sfA_1$).
$$
\epsfbox{cut.4}\qquad\qquad \mapsto \qquad \R\oplus \sfA_1
$$
 Etc., etc.

\SS

4)      The prism has two planes of symmetry. This is by chance,
partially this is induced by
 a symmetry of the initial Coxeter
scheme \ref{dynkin-hyp}. The latter symmetry
implies the equality of {\it strata}:
$$
AB=DE,\quad AC=CE, \quad AD=BC,
\quad BC=CD.
$$

5)      Our prism generates  a reflection group in $\bL^3$.
The reader can easily  imagine a neighborhood of our
prism in $\bL^3$. For instance, near the vertex $A$ we
have the picture drawn in figure \ref{ikos.1}

\SS

6) The development of the prism is a regular
10-gon having right angles; the reflection of
the ``billiard trajectory''
$ABEDA$ is of type d) on figure \ref{billiard}.
 The regularity property follows
by reduction from $\bL^4$, but it is not self-obvious from
the picture of the  3-dimensional prism. Obviously,
diagonals%
\footnote{There are two diagonals $AB$.}
 $AB$ are orthogonal to diagonals $DE$ at
the points of intersection (see the left side of the figure;
but this is not a self-obvious property
of this regular 10-gon).

\SS

7) We observe the second copy of the polygonal line $ADEBA$
in the development. Bending the 10-gon by 
this line, we
obtain a prism congruent to our prism.

In fact, our 10-gon is the picture on the intersection of two
mirrors, denote them by $Y_1$, $Y_2$.
We can roll the simplex $\Sigma$ along each mirror
$Y_1$, $Y_2$ and then we roll it again over the intersection
$Y_1\cap Y_2$. We obtain two different pictures on the 10-gon
and  both  are present in the figure \ref{andreev}.

\nocite{*}\bibliographystyle{amsplain}

\end{document}